\documentclass[11pt]{article}
\usepackage{amssymb}
\newcommand{\be}{\begin{equation}}
\newcommand{\ee}{\end{equation}}
\newcommand{\bea}{\begin{eqnarray}}
\newcommand{\eea}{\end{eqnarray}}

\begin{document}


\title{On the integer partitions recursive structure}
\author{Boris Y. Rubinstein
\\Stowers Institute for Medical Research
\\1000 50th St., Kansas City, MO 64110, U.S.A.}
\date{\today}

\maketitle
\begin{abstract}
Sylvester showed that the partition
of an integer 
into a set of positive integers 
can be represented as a sum of the polynomial term 
and quasiperiodic components called the Sylvester waves.
The wave itself is a weighted sum of the 
polynomial terms multiplied by the periodic functions.
The integer weights are found to be a sum of partitions into a smaller set of integers
implying 
the recursive structure of integer partitions.
\end{abstract}

{\bf Keywords}: integer partitions, Sylvester waves.

{\bf 2010 Mathematics Subject Classification}: 11P82.
\vskip0.4cm
The problem of integer partition into positive integers has a long history started from
the work of Euler who by introducing the idea of generating functions laid 
a foundation of the theory of partitions \cite{GAndrews}.
A. Cayley  and J.J. Sylvester  provided a new insight and made a
remarkable progress in the field. Cayley described the symmetry properties of partitions \cite{Cayley1855}
while Sylvester found 
the procedure for their computation \cite{Sylv1,Sylv2}. Define the 
partition function $W(s,{\bf d}^m)$ is a
number of partitions of integer $s$ into positive integers  ${\bf d}^m=\{d_1,d_2,\ldots,d_m\}$.
Sylvester showed  \cite{Sylv1} that $W(s,{\bf d}^m)$ can be split into polynomial and quasiperiodic parts and presented it as a sum of components
called the {\em Sylvester waves}
\be
W(s,{\bf d}^m) = \sum_{j=1} W_j(s,{\bf d}^m),
\label{SylvWavesExpand}
\ee
where summation runs over all divisors
of the elements in the set ${\bf d}^m$.
The wave $W_j(s,{\bf d}^m)$ is a quasipolynomial in $s$
closely related to prime roots $\rho_j$ of unity.
Namely, Sylvester showed \cite{Sylv2} that the wave
$W_j(s,{\bf d}^m)$ is a coefficient of
${t}^{-1}$ in the series expansion in ascending powers of $t$ of
\be
F_j(s,t)=\sum_{\rho_j} \frac{\rho_j^{-s} e^{st}}{\prod_{k=1}^{m}
        \left(1-\rho_j^{d_k} e^{-d_k t}\right)}\;,
\quad
\rho_j=\exp(2\pi i n/j),
\label{generatorWj}
\ee
where the summation is made over all prime roots of unity
$\rho_j$ for $n$ relatively prime to $j$
(including unity) and smaller than $j$.
Guided by the Sylvester recipe (\ref{generatorWj}) we found in \cite{Rub04} a formula for the
Sylvester wave $W_j(s,{\bf d}^m)$ as a finite sum of the Bernoulli
polynomials of higher
order \cite{bat53,Norlund1924} multiplied by a periodic function of integer 
period $j$. 
The polynomial part $W_1(s,{\bf d}^m)$ of the partition function reads
\begin{equation}
W_1(s,{\bf d}^m) =
\frac{1}{(m-1)!\;\pi_m}
B_{m-1}^{(m)}(s + s_m, {\bf d}^m)\;,
\quad
s_m = \sum_{i=1}^m d_i, \ \ \pi_m =  \prod_{i=1}^m d_i\;,
\label{W_1}
\end{equation}
with the Bernoulli polynomials of higher 
order defined by the generating function \cite{Norlund1924}:
\be
\frac{e^{st} t^m \pi_m}{\prod_{i=1}^m (e^{d_it}-1)} =
\sum_{n=0}^{\infty} B^{(m)}_n(s,{\bf d}^m)
\frac{t^{n}}{n!}\;.
\label{genfuncBernoulli0}
\ee
The Sylvester wave $W_j(s,{\bf d}^m)$ for $j>1$ can be written as
\bea
W_j(s,{\bf d}^m)  =
\frac{j^{k_j-m}}{(m-1)! \; \pi_{m}}
\sum_{{\bf r}=0}^{j-1}
B^{(m)}_{m-1}(s+s_m+ {\bf r}\cdot{\bf d}^{m-k_j},{\bf d}^{m}_j)
\Psi_j(s+s_m+{\bf r}\cdot{\bf d}^{m-k_j}),
\label{WjFin} 
\eea
where we introduce the modified set 
${\bf d}^{m}_j = {\bf d}^{k_j} \cup j{\bf d}^{m-k_j}$
made of two subsets -- 
${\bf d}^{k_j}=\{d_{1},d_{2},\ldots,d_{k_j}\}$ of the generators $d_i$ divisible by $j$
and 
$j{\bf d}^{m-k_j}=\{jd_{k_j+1},jd_{k_j+2},\ldots,jd_{m}\}$ made of $d_i$ nondivisible by $j$.
The $(m-k_j)$-dimensional vector ${\bf r}$
has the form ${\bf r} = \{r_{k_j+1},r_{k_j+2},\ldots,r_{m}\}$
with $0 \le r_i \le j-1$.
The prime circulator $\Psi_j$ in (\ref{WjFin}) is the $j$-periodic function
introduced in \cite{Cayley1855}
that can be written as a sum of simple periodic functions
\be
\Psi_j(s) = 
\sum_{\rho_j} \rho_j^s  = \sum_{n} \psi_{j}^n(s),
\quad
\psi_{j}(s) = \exp(2\pi i s/j),
\quad
\sum_{k=0}^{j-1} \psi_{j}^n(k) = 0,
\label{gencirc}
\ee
for all $1 \le n \le j-1$ coprime to $j$. 
Comparing (\ref{WjFin}) to (\ref{W_1}) we find \cite{Rub04} 
that the Sylvester wave $W_j(s,{\bf d}^m)$ can be viewed as a weighted sum of the 
polynomial parts with shifted argument with the weights given by 
the prime circulator 
\bea
W_j(s,{\bf d}^m)  &=&
\sum_{{\bf r}=0}^{j-1} 
W_1(s-{\bf r}\cdot{\bf d}^{m-k_j},{\bf d}_j^m)
\Psi_j(s-{\bf r}\cdot{\bf d}^{m-k_j})
\label{WjFin0} \\
&=& \sum_{l=0}^{l_{max}} A_l
W_1(s-l,{\bf d}_j^m)
\Psi_j(s-l),
\quad
l_{max} = (j-1) \!\!\!\sum_{i=k_j+1}^{m}\!\! d_i,
\label{WjFin1}\\
l &=& \sum_{{\bf r}=0}^{j-1} {\bf r}\cdot{\bf d}^{m-k_j}
= \sum_{i=k_j+1}^{m} r_i d_i,
\quad
0 \le r_i \le j-1.
\label{main_eq}
\eea
The representation (\ref{WjFin1}) equivalent to (\ref{WjFin0}) 
has integer nonnegative coefficents $A_l$.

Noting that the same 
value $l$ of the sum in (\ref{main_eq})
can be obtained for multiple integer vectors ${\bf r}$
we define $A_l$ in (\ref{WjFin1}) as the number of 
integer vectors ${\bf r}$ satisfying
the Diophantine equation (\ref{main_eq}) together with $m-k_j$ inequalities
$0 \le r_i \le j-1$.
It can be shown (see for example \cite{RubSylvCay,Rub2023}) that
such a problem is equaivalent to the computation of some
integer {\it vector} partition. Namely, introduce $m-k_j$ slack integer variables
$\rho_i,\ k_j+1 \le i \le m$ and 
add to the linear equation in (\ref{main_eq}) $m-k_j$ equations $r_i+\rho_i = j-1$
to form a system of $(m-k_j+1)$
Diophantine equations
\bea
l &=& d_{k_j+1}r_{k_j+1} + d_{k_j+2}r_{k_j+2} +\ldots+ d_{m}r_{m} 
,
\nonumber \\
j-1 &=&  r_{k_j+1} + \rho_{k_j+1},
\nonumber \\
j-1 &=&  r_{k_j+2} + \rho_{k_j+2},
\label{syst}\\
&& \ldots
\nonumber \\
j-1 &=&  r_{m} + \rho_{m}.
\nonumber 
\eea
The equation (\ref{main_eq}) has the same number of solutions
as the system (\ref{syst}) equivalent to the matrix equation
\be
{\bf S} = {\bf D}\cdot{\bf R},
\quad
{\bf D} = 
\left [
\begin{array}{cccccccc}
d_{k_j+1} & d_{k_j+2} & \ldots & d_{m} & 0 & 0 & \ldots & 0 \\
1 & 0 & \ldots & 0 & 1 & 0 & \ldots & 0 \\
0 & 1 & \ldots & 0 & 0 & 1 & \ldots & 0 \\
\ldots & \ldots & \ldots & \ldots & \ldots & \ldots & \ldots & \ldots \\
0 & 0 & \ldots & 1 & 0 & 0 & \ldots & 1 
\end{array}
\right ],
\quad
{\bf S} = 
\left [
\begin{array}{c}
l \\
j-1 \\
j-1 \\
\ldots \\
j-1 
\end{array}
\right ],
\label{matrix}
\ee
where ${\bf D}$ is the  $(m-k_j+1) \times 2(m-k_j)$ integer nonnegative matrix 
while ${\bf S}$ 
and ${\bf R} = \{r_{k_j+1},r_{k_j+2},\ldots,r_{m},\rho_{k_j+1},\rho_{k_j+2},\ldots,\rho_{m}\}^T$
are integer nonnegative vectors.

The idea of the recursive reduction of a vector partition to scalar ones belongs to 
Sylvester \cite{Sylv3} and its first implementation in the case of 
double partitions was performed 
by Cayley \cite{Cayley1860}. 
Sylvester noted in \cite{Sylv3} that this algorithm based on the repeated variable elimination
(corresponding to the elimination of the rows and columns of the matrix ${\bf D}$)
works only when some
additional conditions imposed on the matrix elements
are satisfied. That observation confirmed by Cayley in \cite{Cayley1860}
led mathematicians to a belief of very limited usefulness of this approach and it was eventually forgotten. 
Recently the author of this manuscript generalized the original Sylvester-Cayley method
by lifting all the restrictions and showed that any vector partition can indeed be 
written as a superposition of a finite number 
of scalar partitions \cite{Rub2025}.

Using the approaches presented in \cite{Rub2025} we express
the number $A_l$ of integer nonnegative solutions of the system (\ref{syst}) 
as a sum of the scalar partitions
\bea
A_l = \!\!\!\sum_{p=0}^{2^{m-k_j}-1}\!\!\! (-1)^{\lambda_p}W(s_p,{\bf d}^{m-k_j}),
\quad
\lambda_p = \sum_{i=k_j+1}^{m} v_{p,i},
\quad
s_p = l - j \;{\bf v}_p \cdot {\bf d}^{m-k_j},
\label{A_l}
\eea
where the $(m-k_j)$-dimensional vector ${\bf v}_p=\{v_{p,k_j+1},v_{p,k_j+2},\ldots,v_{p,m}\}$ is 
a sequence of zeros and ones defining the binary representation of
the integer $p$.

This result together with (\ref{WjFin1}) signify 
the recursive nature of the partition function 
structure.
Specifically, the partition function 
$W(s,{\bf d}^m)$ for a $m$-dimensional set of integers ${\bf d}^m$
is a sum of the Sylvester waves $W_j(s,{\bf d}^m)$ for every divisor $j$ of the set ${\bf d}^m$.
The wave $W_j(s,{\bf d}^m)$
is a weighted sum of the $j$-periodically modulated polynomial part $W_1(s-l,{\bf d}_j^m)\Psi_j(s-l)$ with shifted argument.
The weights $A_l$ in their turn are completely determined by the partitions 
$W(s_p,{\bf d}^{m-k_j})$ of $s_p$ into the 
subset  ${\bf d}^{m-k_j}$
containing only the generators nondivisible by $j$. 
Thus the computation of the scalar partition can be viewed as a recursive
process involving partitions with the generator sets of continuously reducing length 
and the problem of integer partitions appears to be {\it self-contained} one.





\end{document}